\documentclass[a4paper,fleqn,10pt]{article}

\usepackage{bbold, epsfig}
\usepackage{latexsym, amsmath, euscript,enumitem}
\usepackage[T1]{fontenc}
\usepackage{fancyhdr}

\newtheorem{theo}{Theorem}

\newcommand{\ds}{\displaystyle}
\renewcommand{\div}{{\rm div \,}}
\newcommand{\grad}{{\rm grad \,}}

\newcommand{\esssup}[1]{\mathop{\rm ess\ sup}}
\newcommand{\essinf}[1]{\mathop{\rm ess\ inf}}
\newcommand{\N}{{\rm I\kern - 2.5pt N}}
\newcommand{\Z}{{\rm Z\kern - 5.5pt Z}}
\newcommand{\Q}{{\rm I\kern - 5.25pt Q}}
\newcommand{\C}{{\rm I\kern - 6.25pt C}}
\newcommand{\R}{{\rm I\kern - 2.5pt R}}
\newcommand{\D}{\mbox \DH}
\begin{document}
\title{On the Maxwell-Stefan approach to multicomponent diffusion}
\author{Dieter Bothe}
\date{}
\maketitle
$\mbox{ }$\vspace{-0.6in}\\
\begin{center}
Center of Smart Interfaces,
TU Darmstadt\\
Petersenstrasse 32,
D-64287 Darmstadt,
Germany\\
e-mail: bothe@csi.tu-darmstadt.de\\[2ex]
\textsl{
Dedicated to Herbert Amann on the occasion of his 70$^{th}\!$ anniversary}
\vspace{0.15in}
\end{center}
\begin{abstract}
\noindent
We consider the system of Maxwell-Stefan equations which describe multicomponent diffusive
fluxes in non-dilute solutions or gas mixtures. We apply the Perron-Frobenius theorem to
the irreducible and quasi-positive matrix which governs the flux-force relations and are
able to show normal ellipticity of the associated multicomponent diffusion operator.
This provides local-in-time wellposedness of the Maxwell-Stefan  multicomponent diffusion
system in the isobaric, isothermal case.
\end{abstract}



\section{Introduction}
On the macroscopic level of continuum mechanical modeling, fluxes of chemical components
(species) are due to convection and molecular fluxes, where the latter essentially refers
to diffusive transport.
The almost exclusively employed constitutive ''law'' to model diffusive fluxes within continuum
mechanical models is Fick's law, stating that the flux of a chemical component is proportional
to the gradient of the concentration of this species, directed against the gradient.
There is no influence of the other components, i.e.\ cross-effects are ignored although
well-known to appear in reality. Actually, such cross-effects can completely divert the diffusive fluxes,
leading to so-called reverse diffusion (up-hill diffusion in direction of the gradient)
or osmotic diffusion (diffusion without a gradient). This has been proven in several experiments,
e.g.\ in a classical setting by  Duncan and Toor; see \cite{Toor62}.

To account for such important phenomena, a multicomponent diffusion approach is required
for realistic models. The standard approach within the theory of
Irreversible Thermodynamics replaces Fickian fluxes by linear combinations of the gradients
of all involved concentrations, respectively chemical potentials. This requires the knowledge of a full matrix of binary diffusion
coefficients and this diffusivity matrix has to fulfill certain requirements like positive
semi-definiteness in order to be consistent with the fundamental laws from thermodynamics.
The Maxwell-Stefan approach to multicomponent diffusion leads to a concrete form of the diffusivity
matrix and is based on molecular force balances to relate all individual species velocities.
While the Maxwell-Stefan equations are successfully used in engineering applications, they seem
much less known in the mathematical literature.
In fact we are not aware of a rigorous mathematical analysis of the Maxwell-Stefan approach
to multicomponent diffusion, except for \cite{Giovan} which mainly addresses questions of modeling and numerical computations, but also contains some analytical results which are closely related to the present considerations.

\section{Continuum Mechanical Modeling\\
of Multicomponent Fluids}
%\noindent
We consider a multicomponent fluid composed of $n$ chemical components $A_i$.
Starting point of the Maxwell-Stefan equations are the individual mass balances, i.e.
\begin{equation}
\label{mass-bal}
\partial_t \rho_i + \div (\rho_i {\bf u}_i)=R_i^{\rm tot},
\end{equation}
where $\rho_i=\rho_i (t,{\bf y})$ denotes the mass density and ${\bf u}_i={\bf u}_i (t,{\bf y})$ the individual velocity of species $A_i$.
Note that the spatial variable is denoted as ${\bf y}$, while the usual symbol ${\bf x}$ will refer to the composition of the mixture.
The right-hand side is the total rate of change of species mass due to all chemical transformations.
We assume conservation of the total mass, i.e.\ the production terms satisfy $\sum_{i=1}^n R_i^{\rm tot} =0$.
Let $\rho$ denote the total mass density and ${\bf u}$ be the barycentric
(i.e., mass averaged) velocity, determined by
\[
\rho :=\sum_{i=1}^n \rho_i, \qquad \rho \, {\bf u}:=\sum_{i=1}^n \rho_i \, {\bf u}_i.
\]
Summation of the individual mass balances
\eqref{mass-bal} then yields
\begin{equation}
\label{conti-eq}
\partial_t \rho + \div (\rho {\bf u})=0,
\end{equation}
i.e.\ the usual continuity equation.

In principle, a full set of $n$ individual momentum balances should now be added to the model; cf.\ \cite{Kerk}.
But in almost all engineering models, a single set of
Navier-Stokes equations is used to describe the evolution of the velocity field, usually without
accounting for individual contributions to the stress tensor.
One main reason is a lack of information about appropriate constitutive equations for the stress
in multicomponent mixtures; but cf.\ \cite{Raja}.
For the multicomponent, single momentum model the barycentric velocity
${\bf u}$ is assumed to be determined by the Navier-Stokes equations.
Introducing the mass diffusion fluxes
\begin{equation}
\label{mass-flux}
{\bf j}_i:=\rho_i ({\bf u}_i -{\bf u})
\end{equation}
and the mass fractions $Y_i :=\rho_i /\rho$, the mass balances \eqref{mass-bal} can be rewritten as
\begin{equation}
\label{mass-frac-bal}
\rho \partial_t Y_i + \rho  {\bf u} \cdot \nabla Y_i  + \div {\bf j}_i=R_i^{\rm tot}.
\end{equation}
In the present paper, main emphasis is on the aspect of multicomponent diffusion, including the cross-diffusion effects.
Therefore, we focus on the special case of \emph{isobaric, isothermal} diffusion. The (thermodynamic) pressure $p$
is the sum of partial pressures $p_i$ and the latter correspond to $c_i R\, T$ in the general case with $c_i$
denoting the molar concentration, $R$ the universal gas constant and $T$ the absolute temperature;
here $c_i=\rho_i /M_i$ with $M_i$ the molar mass of species $A_i$.
Hence isobaric conditions correspond to the case of constant total molar concentration
 $c_{\rm tot}$, where $c_{\rm tot}:=\sum_{i=1}^n c_i$.
Still, species diffusion can lead to transport of momentum because the $M_i$ are different.
Instead of ${\bf u}$ we therefore employ the molar averaged velocity defined by
\begin{equation}
\label{velo-def}
c_{\rm tot} {\bf v}:=\sum_{i=1}^n c_i {\bf u}_i.
\end{equation}
Note that other velocities are used as well; only the diffusive fluxes have to be adapted;
see, e.g., \cite{TK-book}.
With the molar averaged velocity, the species equations (\ref{mass-bal}) become
\begin{equation}
\label{species-bal}
\partial_t c_i + \div (c_i {\bf v} + {\bf J}_i)=r_i^{\rm tot}
\end{equation}
with $r_i^{\rm tot}:=R_i^{\rm tot}/M_i$ and the diffusive molar fluxes
\begin{equation}
\label{mol-flux}
{\bf J}_i:=c_i ({\bf u}_i -{\bf v}).
\end{equation}
Below we exploit the important fact that
\begin{equation}
\label{flux-sum}
\sum_{i=1}^n {\bf J}_i =0.
\end{equation}
As explained above we may now assume ${\bf v}=0$ in the isobaric case.
In this case the species equations \eqref{species-bal} simplify to a system of
reaction-diffusion systems given by
\begin{equation}
\label{species-bal-diff}
\partial_t c_i + \div {\bf J}_i=r_i^{\rm tot},
\end{equation}
where the individual fluxes ${\bf J}_i$ need to be modeled by appropriate constitutive
equations.
The most common constitutive equation is Fick's law which states that
\begin{equation}
\label{Fick}
{\bf J}_i = - D_i \, \grad c_i
\end{equation}
with diffusivities $D_i >0$. The diffusivities are usually assumed to be constant,
while they indeed depend in particular on the composition of the system, i.e.\ $D_i = D_i({\bf c})$ with
${\bf c}:=(c_1\ldots ,c_n)$. Even if the dependence of the $D_i$ is taken into account, the above
definition of the fluxes misses the cross-effects between the diffusing species.
In case of concentrated systems more realistic constitutive equations are hence required which
especially account for such mutual influences.
Here a common approach is the general constitutive law
\begin{equation}
\label{MCD}
{\bf J}_i = - \sum_{j=1}^n D_{ij} \, \grad c_j
\end{equation}
with binary diffusivities $D_{ij}=D_{ij}({\bf c})$.
Due to the structure of the driving forces, as discussed below, the matrix
${\bf D}=[D_{ij}]$ is of the form ${\bf D}({\bf c}) = {\bf L}({\bf c}) \, G''({\bf c})$
with a positive definite matrix $G''({\bf c})$, the Hessian of the Gibbs free energy.
Then, from general principles
of the theory of Irreversible Thermodynamics,
it is assumed that the matrix of transport coefficients ${\bf L}=[L_{ij}]$ satisfies
\begin{itemize}
\item ${\bf L}$ is {\it symmetric} (the Onsager reciprocal relations)
\item ${\bf L}$ is {\it positive semidefinite} (the second law of thermodynamics).
\end{itemize}
Under this assumption the quasilinear reaction-diffusion system
\begin{equation}
\label{species-bal-sys}
\partial_t {\bf c} + \div (-{\bf D}({\bf c})\, \nabla {\bf c})={\bf r}({\bf c}),
\end{equation}
satisfies - probably after a reduction to $n-1$ species - parabolicity conditions sufficient for local-in-time wellposedness.
Here ${\bf r}({\bf c})$ is short for $(r_1^{\rm tot} ({\bf c}), \ldots ,r_n^{\rm tot} ({\bf c}))$.

A main problem now is how realistic diffusivity matrices together with their dependence on the composition vector ${\bf c}$ can be obtained.

Let us note in passing that Herbert Amann has often been advocating that general flux vectors should be considered, accounting both for concentration dependent diffusivities and for cross-diffusion effects.
For a sample of his contributions to the theory of reaction-diffusion systems with general flux vectors
see \cite{Amann-RD},  \cite{Amann1} and the references given there.

\section{The Maxwell-Stefan Equations}
The Maxwell-Stefan equations rely on inter-species force balances.
More precisely, it is assumed that the thermodynamical driving force
${\bf d}_i$ of species $A_i$ is in local equilibrium with the total friction force.
Here and below it is often convenient to work with the molar fractions
$x_i:=c_i / c_{\rm tot}$ instead of the chemical concentrations.
From chemical thermodynamics it follows that for multicomponent systems which are locally close to thermodynamical
equilibrium (see, e.g., \cite{TK-book}) the driving forces under isothermal conditions are given as
\begin{equation}
\label{driving}
{\bf d}_i = \frac{x_i}{R\, T}\grad \mu_i
\end{equation}
with $\mu_i$ the chemical potential
of species $A_i$. Equation \eqref{driving} requires some more explanation.
Recall first that the chemical potential $\mu_i$ for species $A_i$ is defined as
\begin{equation}
\label{chem-pot}
\mu_i = \frac{\partial G}{\partial c_i},
\end{equation}
where $G$ denotes the
(volume-specific) density of the Gibbs free energy. The chemical potential depends on $c_i$,
but also on all other $c_j$ as well as on pressure  and temperature.
In the engineering literature, from the chemical potential a part $\mu_i^0$ depending on pressure and temperature
is often separated and, depending on the context, a gradient may be applied only to the remainder.
To avoid confusion, the common notation in use therefore is
\[
\nabla \mu_i =\nabla_{T,p} \mu_i + \frac{\partial \mu_i}{\partial p} \nabla p + \frac{\partial \mu_i}{\partial T} \nabla T.
\]
Here $\nabla_{T,p} \mu_i$ means the gradient taken under constant pressure and temperature.
In the isobaric, isothermal case this evidently makes no difference.
Let us also note that $G$ is assumed to be a convex function of the $c_i$ for
single phase systems, since this guarantees thermodynamic stability, i.e.\ no spontaneous
phase separations.
For concrete mixtures, the chemical potential is often assumed to be given by
\begin{equation}
\label{chempot}
\mu_i = \mu_i^0 + R\, T \ln a_i
\end{equation}
with $a_i$ the so-called activity of the $i$-th species; equation (\ref{chempot})
actually implicitly defines $a_i$. In (\ref{chempot}), the term $\mu_i^0$ depends on pressure and temperature.
For a mixture of ideal gases, the activity $a_i$ equals
the molar fraction $x_i$. The same holds for solutions in the limit of an ideally dilute component, i.e.\ for $x_i\to 0+$.
This is no longer true for non-ideal systems in which case the activity is written as
\begin{equation}
\label{activity}
a_i = \gamma_i \, x_i
\end{equation}
with an activity coefficient $\gamma_i$ which itself depends in particular on the full
composition vector ${\bf x}$.

The mutual friction force between species $i$ and $j$ is assumed to be
proportional to the relative velocity as well as to the amount of molar mass.
Together with the assumption of balance of forces this leads to the relation
\begin{equation}
\label{force-balance}
{\bf d}_i = - \sum_{j\neq i} f_{ij} \, x_i \, x_j ({\bf u}_i -{\bf u}_j)
\end{equation}
with certain drag coefficients $f_{ij}>0$; here $f_{ij}=f_{ji}$ is a natural mechanical assumption.
Insertion of (\ref{driving}) and introduction of
the so-called Maxwell-Stefan (MS) diffusivities $\D_{ij}=1/f_{ij}$ yields the system
\begin{equation}
\label{MS-system}
 \frac{x_i}{R\, T}\grad \mu_i = -\, \sum_{j\neq i} \frac{ x_j {\bf J}_i - x_i {\bf J}_j}{c_{\rm tot}\, \D_{ij}}
\quad \mbox{ for } i=1,\ldots ,n.
\end{equation}
The set of equations (\ref{MS-system}) together with (\ref{flux-sum})
forms the Maxwell-Stefan equations of multicomponent diffusion.
The matrix $[\D_{ij}]$ of MS-diffusi\-vities is assumed to be symmetric in accordance
with the symmetry of $[f_{ij}]$.
Let us note that for ideal gases the symmetry can be obtained from the kinetic theory
of gases; cf.\ \cite{Hirsch} and \cite{Muck}. The MS-diffusi\-vities $\D_{ij}$ will in general depend on the
composition of the system.

Due to the symmetry of $[\D_{ij}]$, the model is in fact consistent with the Onsager reciprocal relations
(cf.\ \cite{Standart} as well as below), but notice that the
$\D_{ij}$ are not to be inserted into \eqref{MCD}, i.e.\ they do not directly correspond to
the $D_{ij}$ there. Instead, the MS equations have to be inverted in order to provide the
fluxes ${\bf J}_i$.

Note also that the Ansatz \eqref{force-balance} implies $\sum_i {\bf d}_i =0$
because of the symmetry of $[f_{ij}]$, resp.\ of $[\D_{ij}]$.
Hence $\sum_i {\bf d}_i =0$ is necessary in order for \eqref{force-balance} to be consistent.
It in fact holds because of (and is nothing but) the Gibbs-Duhem relation, see e.g.\ \cite{KW}.
The relation $\sum_i {\bf d}_i =0$ will be important below.
\vspace{0.1in}\\
{\bf Example (Binary systems).} For a system with two components we have
\begin{equation}
\label{binary1}
{\bf d}_1 (=- {\bf d}_2 )= - \, \frac{1}{c_{\rm tot} \D_{12}}\big( x_2 {\bf J}_1 - x_1 {\bf J}_2 \big).
\end{equation}
Using $x_1 + x_2=1$ and ${\bf J}_1+{\bf J}_2=0$ one obtains
\begin{equation}
\label{binary2}
{\bf J}_1 (=- {\bf J}_2 )= -\frac{\D_{12}}{R \, T}c_1 \, \grad \mu_1.
\end{equation}
Writing $c$ and ${\bf J}$ instead of $c_1$ and ${\bf J}_1$, respectively, and assuming that the
chemical potential is of the form $\mu = \mu^0 + R\, T \ln (\gamma c)$ with
the activity coefficient $\gamma = \gamma (c)$ this finally yields
\begin{equation}
\label{binary3}
{\bf J}= -\D_{12}  \Big( 1+\frac{c \, \gamma' (c)}{\gamma (c)} \Big) \grad c .
\end{equation}
Inserting this into the species equation leads to a nonlinear diffusion equation,
namely
\begin{equation}
\label{filtration}
\partial_t c - \Delta \phi (c) =r(c),
\end{equation}
where the function $\phi :\R \to \R$ satisfies $\phi' (s)=\D_{12}(1+s \gamma' (s)/\gamma (s))$ and, say, $\phi (0)=0$.
Equation (\ref{filtration}) is also known as the filtration equation (or, the generalized
porous medium equation) in other applications.
Note that well-known pde-theory applies to \eqref{filtration} and especially provides
well-posedness as soon as $\phi$ is continuous and nondecreasing; cf., e.g., \cite{Vazquez}.
The latter holds if $s\to s\gamma(s)$ is increasing
which is nothing but the fact that the chemical potential $\mu$ of a component should be an increasing function of its concentration. This is physically reasonable in systems without
phase separation.

\section{Inversion of the Flux-Force Relations}
\label{MS-inversion}
In order to get constitutive equations for the fluxes ${\bf J}_i$
from the Maxwell-Stefan equations,
which need to be inserted into \eqref{species-bal-diff}, we have to invert \eqref{MS-system}.
Now \eqref{MS-system} alone is not invertible for the fluxes, since these are linearly dependent.
Elimination of ${\bf J}_n$ by means of (\ref{flux-sum}) leads to the reduced system
\begin{equation}
\label{reduced-system}
c_{\rm tot}
\left[
\begin{array}{ccc}
\; & {\bf d}_1 & \; \\[-1ex]
& \cdot &\\[-2ex]
& \cdot &\\[-2ex]
& \cdot &\\[-1ex]
\; & {\bf d}_{n-1} & \;
\end{array}
\right]
= -\,{\bf B}
\left[
\begin{array}{ccc}
\; & {\bf J}_1 & \; \\[-1ex]
& \cdot &\\[-2ex]
& \cdot &\\[-2ex]
& \cdot &\\[-1ex]
\; & {\bf J}_{n-1} & \;
\end{array}
\right],
\end{equation}
where the $(n-1)\times (n-1)$-matrix ${\bf B}$ is given by
\begin{equation}
\label{B-def}
B_{ij}=
\left \{
\begin{array}{ll}
\ds
x_i\Big(\frac{1}{\D_{1n}} - \frac{1}{\D_{ij}} \Big) & \mbox{for } i\neq j,\\[3ex]
\ds
\frac{x_i}{\D_{in}} + \sum_{k\neq i}^n \frac{x_k}{\D_{ik}}& \mbox{for } i= j
\; \mbox{ (with $x_n=1-\sum_{m<n} x_m$)}.
\end{array}
\right .
\end{equation}
Assuming for the moment the invertibility of ${\bf B}$ and letting $\mu_i$ be functions of
the composition expressed by the molar fractions ${\bf x}=(x_1, \ldots ,x_n)$, the fluxes are
given by
\begin{equation}
\label{MS-fluxes}
\left[
\begin{array}{ccc}
\; & {\bf J}_1 & \; \\
& \cdot &\\
& \cdot &\\
& \cdot &\\
\; & {\bf J}_{n-1} & \;
\end{array}
\right]
=
- c_{\rm tot} {\bf B}^{-1}\, {\bf \Gamma} \,
\left[
\begin{array}{ccc}
\; & \nabla x_1 & \; \\
& \cdot &\\
& \cdot &\\
& \cdot &\\
\; & \nabla x_{n-1} & \;
\end{array}
\right],
\end{equation}
where
\begin{equation}
\label{G-def}
{\bf \Gamma} = [\Gamma_{ij}]
\quad \mbox{ with } \;
\Gamma_{ij} = \delta_{ij} + x_i \, \frac{\partial \ln \gamma_i}{\partial x_j}
\end{equation}
captures the thermodynamical deviations from the ideally diluted situation;
here $\delta_{ij}$ denotes the Kronecker symbol.
\vspace{0.1in}\\
{\bf Example (Ternary systems).} We have
\begin{equation}
\label{B-ternary}
{\bf B}=
\left[
\begin{array}{cc}
\frac{1}{\D_{13}}+x_2\Big(\frac{1}{\D_{12}} -\frac{1}{\D_{13}}\Big)
  & - x_1\Big(\frac{1}{\D_{12}} -\frac{1}{\D_{13}}\Big) \\[2ex]
- x_2\Big(\frac{1}{\D_{12}} -\frac{1}{\D_{23}}\Big)
 & \frac{1}{\D_{23}}+x_1\Big(\frac{1}{\D_{12}} -\frac{1}{\D_{23}}\Big)
\end{array}
\right]
\end{equation}
and $\det ({\bf B}-t {\bf I})=t^2 - {\rm tr\, } {\bf B} \, t + \det {\bf B}$ with
\begin{equation}
\label{detB-ternary}
{\rm det}\, {\bf B}=
\frac{x_1}{\D_{12}\, \D_{13}}+\frac{x_2}{\D_{12}\, \D_{23}}+\frac{x_3}{\D_{13}\, \D_{23}}
\geq \min \{\frac{1}{\D_{12}\, \D_{13}},\frac{1}{\D_{12}\, \D_{23}},\frac{1}{\D_{13}\, \D_{23}}\}
\end{equation}
and
\begin{equation}
\label{coeff-ternary}
{\rm tr\, } {\bf B} = \frac{x_1+x_2}{\D_{12}}+\frac{x_1+x_3}{\D_{13}}+\frac{x_2+x_3}{\D_{23}} \geq
2 \min \{\frac{1}{\D_{12}},\frac{1}{\D_{13}},\frac{1}{\D_{23}}\}.
\end{equation}
It is easy to check that $({\rm tr\, } {\bf B})^2 \geq 3 \det {\bf B}$ for this
particular matrix and therefore the spectrum of ${\bf B}^{-1}$ is in the right complex half-plane
within a sector of angle less than $\pi /6$. This implies normal ellipticity of the
differential operator ${\bf B}^{-1}({\bf x}) (-\Delta {\bf x})$.
Recall that a second order differential operator with matrix-valued coefficients
is said to be \emph{normally elliptic} if the symbol of the principal part has it's spectrum inside the open right
half-plane of the complex plane; see section~4 in \cite{Amann1} for more details. This notion has been introduced by
Herbert Amann in \cite{Amann-RD}
as the appropriate concept for generalizations to more general situations with operator-valued coefficients.

Consequently, the Maxwell-Stefan equations for a ternary system are locally-in-time wellposed
if ${\bf \Gamma}={\bf I}$, i.e.\ in the special case of \emph{ideal solutions}.
The latter refers to the case when the chemical potentials are of
the form (\ref{chempot}) with $\gamma_i \equiv 1$ for all $i$.
Of course this extends to any ${\bf \Gamma}$ which is a small perturbations of ${\bf I}$, i.e.\
to slightly non-ideal solutions.

Let us note that Theorem~1 below yields the local-in-time wellposedness also for
general non-ideal solutions provided the Gibbs energy is strongly convex.
Note also that the reduction to $n-1$ species is the common approach in the
engineering literature, but invertibility of ${\bf B}$ is not rigorously checked.
For $n=4$, the $3\times 3$-matrix ${\bf B}$ can still be shown to be invertible for any composition
due to $x_i \geq 0$ and $\sum_i x_i =1$. Normal ellipticity can no longer be seen so easily.
For general $n$ this approach is not feasible and the invariant approach below is
preferable.

Valuable references for the Maxwell-Stefan equations and there applications in the
Engineering Sciences are in particular
the books \cite{Bird}, \cite{Hirsch}, \cite{TK-book} and the review article \cite{KW}.

\section{Wellposedness of the Maxwell-Stefan equations}
\label{wellposed}
We first invert the Maxwell-Stefan equations using an invariant formulation.
For this purpose, recall that $\sum_i u_i =0$ holds for both $u_i ={\bf J}_i$ and $u_i ={\bf d}_i$.
We therefore have to solve
\begin{equation}
\label{MS-inv}
 A\, {\bf J} = c_{\rm tot} \, {\bf d}\qquad \mbox{ in } E=\{u\in \R^n : \sum_i u_i =0\},
\end{equation}
$\mbox{ }$\vspace{-0.25in}\\
where $A=A({\bf x})$ is given by
\[
A=
\left[
\begin{array}{ccccc}
-s_1 & \, & \, & \, & d_{ij} \\[-1ex]
\; & \cdot & \, &  & \, \\[-2ex]
\; & \, & \cdot & \, & \;\\[-2ex]
\; &  & \, & \cdot & \; \\[-2ex]
\; d_{ij} & \, & \, & \, & -s_n
\end{array}
\right]
\qquad \mbox{ with }
s_i=\sum_{k\neq i} \frac{x_k}{\D_{ik}},\quad
d_{ij}=\frac{x_i}{\D_{ij}}.
\]
The matrix $A$ has the following properties, where ${\bf x}\gg 0$ means $x_i >0$ for
all $i$:
\begin{enumerate}[label=(\roman{enumi})]
%\begin{enumerate}
\item
$N(A)={\rm span} \{{\bf x}\}$ for ${\bf x} = (x_1,\ldots ,x_n)$.
\item
$R(A)=\{{\bf e}\}^\perp$ for ${\bf e} = (1,\ldots ,1)$.
\item
$A=[a_{ij}]$ is quasi-positive, i.e.\ $a_{ij}\geq 0$ for $i\neq j$.
\item
If ${\bf x}\gg0$ then $A$ is irreducible, i.e.\ for every disjoint partition
$I\cup J$ of $\{1,\ldots ,n\}$ there is some $(i,j)\in I\times J$ such that $a_{ij}\neq 0$.
\end{enumerate}
Due to (i) and (ii) above, the Perron-Frobenius theorem in the version
for quasi-positive matrices applies; cf.\ \cite{Horn} or \cite{Serre}.
This yields the following properties of the spectrum  $\sigma (A)$:
The spectral bound $s(A):=\max \{{\rm Re}\, \lambda : \lambda \in \sigma (A)\}$
is an eigenvalue of $A$, it is in fact a simple eigenvalue with a strictly positive eigenvector.
All other eigenvalues do not have positive eigenvectors or
positive generalized eigenvectors. Moreover,
\[
{\rm Re}\, \lambda < s(A) \quad \mbox{ for all }
\lambda \in \sigma (A), \, \lambda \neq s(A).
\]
From now on we assume that in the present case ${\bf x}$ is strictly positive.
Then, since ${\bf x}$ is an eigenvector to the eigenvalue 0, it follows that
\[
\sigma (A) \subset \{0\}\cup \{z\in \C: {\rm Re}\, z < 0\}.
\]
Unique solvability of \eqref{MS-inv} already follows at this point. In addition, the same
arguments applied to $A_\mu := A- \mu ( {\bf x}\otimes {\bf e} )$ for $\mu \in \R$ yield
\[
\sigma (A_\mu) \subset \{- \mu\}\cup \{z\in \C: {\rm Re}\, z < -\mu \}
\quad \mbox{ for all small } \mu >0.
\]
In particular, $A_\mu$ is invertible for sufficiently small $\mu >0$ and
\begin{equation}
{\bf J} = - c_{\rm tot} \,  \big( A- \mu  ({\bf x}\otimes {\bf e}) \big)^{-1}  {\bf d}
\end{equation}
is the unique solution of \eqref{MS-inv}. Note that $A_\mu {\bf y} ={\bf d}$ with
${\bf d}\perp {\bf e}$ implies ${\bf y}\perp {\bf e}$ and $A {\bf y}= {\bf d}$.
A similar representation of the inverted Maxwell-Stefan
equations can be found in \cite{Giovan}.

The information on the spectrum of $A$ can be significantly improved by symmetrization.
For this purpose let $X={\rm diag} (x_1,\ldots x_n)$ which is regular due to
${\bf x}\gg 0$. Then $A_S:=X^{-\frac 1 2}\,A\, X^{\frac 1 2}$ satisfies
\[
A_S=
\left[
\begin{array}{ccccc}
-s_1 & \, & \, & \, & \hat{d}_{ij} \\[-1ex]
\; & \cdot & \, &  & \, \\[-2ex]
\; & \, & \cdot & \, & \;\\[-2ex]
\; &  & \, & \cdot & \; \\[-2ex]
\; \hat{d}_{ij} & \, & \, & \, & -s_n
\end{array}
\right],
\quad
s_i=\sum_{k\neq i} \frac{x_k}{\D_{ik}},\quad
\hat{d}_{ij}=\frac{\sqrt{x_i x_j}}{\D_{ij}},
\]
i.e.\ $A_S$ is symmetric with $N(A_S)={\rm span} \{\sqrt{\bf x} \}$, where $\sqrt{\bf x}_i:=
\sqrt{x_i}$. Hence the spectrum of $A_S$ and, hence, that of $A$ is real.
Moreover,
\[
A_S(\alpha )=A_S - \alpha \sqrt{\bf x} \otimes \sqrt{\bf x}
\]
has the same properties as $A_S$ for sufficiently small $\alpha >0$.
In particular, $A_S$ is quasi-positive, irreducible and $\sqrt{\bf x} \gg 0$ is an eigenvector
for the eigenvalue $- \alpha$.
This holds for all $\alpha < \delta:=\min \{1/ \D_{ij}: i\neq j\}$.
Hence we obtain the improved inclusion
\begin{equation*}
\sigma (A) \setminus\{0\}=\sigma (A_S (\alpha)) \setminus\{- \alpha \}
\;\mbox{ for all } \alpha \in [0,\delta ).
\end{equation*}
Therefore
\begin{equation}
\label{spectral-incl}
\sigma (A) \subset  (-\infty ,-\delta] \cup \{ 0\},
\end{equation}
which provides a uniform spectral gap for $A$ sufficient to obtain normal ellipticity of
the associated differential operator.

In order to work in a subspace of the composition space $\R^n$ instead of a hyperplane,
let $u_i = c_i - c_{\rm tot}^0 /n$ such that $\sum_i c_i \equiv const$
is the same as $u \in E=\{u\in \R^n : \sum_i u_i =0\}$.
Above we have shown in particular that $A_{|E}:E \to E$ is invertible and
\begin{equation}
\label{flux-repr1}
\left[ {\bf J}_i \right]
=
X^{\frac 1 2} (A_{S|\hat{E}})^{-1} X^{-\frac 1 2}
\left[
{\bf d}_i
\right]
=
\frac{1}{R T} X^{\frac 1 2} (A_{S|\hat{E}})^{-1} X^{\frac 1 2}
\left[
\nabla \mu_i
\right]
\end{equation}
with the symmetrized form $A_S$ of $A$ and $\hat{E}:=X^{\frac 1 2} E=\{\sqrt{\bf x} \}^\perp$.
Note that this also shows the consistency with
the Onsager relations. To proceed, we employ \eqref{chem-pot} to obtain the representation
\begin{equation}
\label{flux-repr}
\left[ {\bf J}_i \right]
=
 X^{\frac 1 2} (A_{S|\hat{E}})^{-1} X^{\frac 1 2}
G''({\bf x}) \, \nabla {\bf x}.
\end{equation}
Inserting \eqref{flux-repr} into \eqref{species-bal-diff} and using $c_{\rm tot} x_i=u_i + c_{\rm tot}^0 /n$,
we obtain the system of species equations with multicomponent diffusion modeled by the Maxwell-Stefan equations.
Without chemical reactions and in an isolated
domain $\Omega \subset \R^n$  (with $\nu$ the outer normal)
we obtain the initial boundary value problem
\begin{equation}
\label{IBVP}
\partial_t u +\div (- {\bf D}(u)\nabla u) =0,
\quad \partial_\nu u_{|\partial \Omega} =0, \; u_{|t=0} =u_0,
\vspace{0.05in}
\end{equation}
which we will consider in $L^p (\Omega; E)$. Note that
$ X^{\frac 1 2} (A_{S|\hat{E}})^{-1} X^{\frac 1 2}
G''({\bf x})$ from \eqref{flux-repr} corresponds to  $- {\bf D}(u)$ here.

Applying well-known results for quasilinear parabolic systems based on $L_p$-maximal regularity,
e.g.\ from \cite{Amann-2005} or \cite{Jan-Bari}, we obtain the following result on local-in-time wellposedness
of the Maxwell-Stefan equations in the isobaric, isothermal case.
Below we call $G\in C^2(V)$ strongly convex if $G''({\bf x})$ is positive definite for all ${\bf x}\in V$.
\begin{theo}
Let $\Omega \subset \R^N$ with $N\geq 1$ be open bounded with smooth $\partial \Omega$.
Let $p>\frac{N+2}{2}$ and $u_0\in W_p^{2-\frac{2}{p}} (\Omega ;E)$ such that $c_i^0 >0$ in $\bar{\Omega}$ and $c_{\rm tot}^0$ is constant in $\Omega$. Let the diffusion matrix ${\bf D}(u)$ be given according to \eqref{flux-repr}, i.e.\ by
\[
{\bf D}(u)= X^{\frac 1 2} (A_{S|\hat{E}})^{-1} X^{\frac 1 2} G''({\bf x}) \;\mbox{ with } c_{\rm tot} x_i=u_i + c_{\rm tot}^0 /n,
\]
where $G:(0,\infty)^n\to \R$ is smooth and strongly convex.
Then there exists - locally in time - a unique strong solution (in the $L^p$-sense) of
\eqref{IBVP}. This solution is in fact classical.
\vspace{0.05in}
\end{theo}
Concerning the proof let us just mention that
\[
\div (- {\bf D}(u)\nabla u)={\bf D}(u)\,(- \Delta u)\; + \mbox{ lower order terms},
\]
hence the system of Maxwell-Stefan equations is locally-in-time wellposed
if the principal part ${\bf D}(u)\,(- \Delta u)$ is normally elliptic for all $u\in E$ such that
${\bf c}(u):=u+c^0_{\rm tot} {\bf e}$ is close to ${\bf c}^0$.
The latter holds if, for some angle $\theta \in (0,\frac \pi 2)$, the spectrum of ${\bf D}(u)\in {\cal L} (E)$ satisfies
\begin{equation}
\label{spectral-cond}
\sigma ({\bf D}(u)) \subset \Sigma_\theta := \{\lambda \in \C\setminus \{0\}: |{\rm arg}\, \lambda |
<\theta\}
\end{equation}
for all $u\in E$ such that $|{\bf c}(u)-{\bf c}^0|_\infty <\epsilon $ for $\epsilon :=\min_i c_i^0 /2$, say.
For such an $u\in E$, let $\lambda \in \C$ and $v\in E$ be such that ${\bf D}(u)\, v = \lambda v$.
Let ${\bf x}:={\bf c}(u)/c_{\rm tot}(u)\in (0,\infty)^n$ and $X={\rm diag}(x_1,\ldots ,x_n)$. Then
\[
X^{\frac 1 2} (A_{S|\hat{E}})^{-1} X^{\frac 1 2} G''({\bf x})\, v= \lambda v.
\]
Taking the inner product with $G''({\bf x})\, v$ yields
\[
\langle (A_{S|\hat{E}})^{-1} X^{\frac 1 2} G''({\bf x})\, v, X^{\frac 1 2} G''({\bf x})\, v\rangle
= \lambda \langle v, G''({\bf x})\, v\rangle.
\]
Note that $X^{\frac 1 2} G''({\bf x})\, v \in \{\sqrt{\bf x} \}^\perp$, hence the left-hand side is strictly positive
due to the analysis given above. Moreover $\langle v, G''({\bf x})\, v\rangle >0$ since $G$ is strongly convex,
hence $\lambda >0$. This implies \eqref{spectral-cond} for any $\theta \in (0,\frac \pi 2)$ and,
hence, local-in-time existence follows.

\section{Final Remarks}
A straight-forward extension of Theorem~1 to the inhomogeneous case
with locally Lipschitz continuous right-hand side $f:\R^n \to \R^n$, say, is possible
if $f(u)\in E$ holds for all $u$. Translated back to the original variables (keeping
the symbol $f$) this yields
a local-in-time solution of
\begin{align*}
\partial_t {\bf c} +\div (- {\bf D}({\bf c})\nabla {\bf c}) = f({\bf c}),
\qquad \partial_\nu {\bf c}_{|\partial \Omega} =0, \quad {\bf c}_{|t=0} ={\bf c}_0
\end{align*}
for appropriate initial values ${\bf c}_0$. Then a natural question is whether the solution
stays componentwise nonnegative. This can only hold if $f$ satisfies
\[
f_i ({\bf c}) \geq 0 \;\mbox{ whenever  $\, {\bf c} \geq 0$ with } c_i=0,
\]
which is called {\it quasi-positivity} as in the linear case.
In fact, under the considered assumption, quasi-positivity of $f$ forces any classical solution to
stay nonnegative as long as it exists. The key point here is the structure of the
Maxwell-Stefan equations \eqref{MS-system} which yields
\begin{align*}
{\bf J}_i =- D_i ({\bf c}) \, \grad c_i + c_i \, {\bf F}_i ({\bf c}, \grad {\bf c})
\end{align*}
with
\begin{align*}
D_i ({\bf c})=1/  \sum_{j\neq i} \frac{ x_j } { \D_{ij}}
\; \mbox{ and }\; {\bf F}_i ({\bf c}, \grad {\bf c})=D_i ({\bf c})  \sum_{j\neq i}\frac{1}{ \D_{ij}} {\bf J}_j.
\end{align*}
Note that $D_i ({\bf c})> 0$ and ${\bf J}_i$ becomes proportional to $\grad c_i$ at points where $c_i$ vanishes, i.e.\
the diffusive cross-effects disappear. Moreover, it is easy to check that
\[
\div {\bf J}_i = D_i ({\bf c}) \Delta c_i \geq 0
\;\mbox{ if $c_i=0$ and } \grad c_i=0.
\]
%\pause
To indicate a rigorous proof for the nonnegativity of solutions, consider the modified system
\begin{equation}
\label{mod-system}
\partial_t c_i +\div {\bf J}_i ({\bf c}) =f_i(t,{\bf c^+}) +\epsilon,
\qquad \partial_\nu {\bf c}_{|\partial \Omega} =0, \quad {\bf c}_{|t=0} ={\bf c}_0 +\epsilon {\bf e},
\end{equation}
where $r^+ :=\max \{r,0\}$ denotes the positive part. Assume that the right-hand side $f$ is quasi-positive and that
\eqref{mod-system} has a classical solution ${\bf c}^\epsilon$ for all small $\epsilon >0$ on a common
time interval $[0,T)$.
Now suppose that, for some $i$, the function
$m_i (t)=\min_{{\bf y}\in \bar{\Omega}} c_i^\epsilon (t,{\bf y})$ has a first zero at $t_0\in (0,T)$.
Let the minimum of $c_i^\epsilon (t_0,\cdot)$ be attained at ${\bf y}_0$ and assume first that ${\bf y}_0$ is an
interior point. Then $c_i^\epsilon (t_0,{\bf y}_0)=0$, $\partial_t c_i^\epsilon (t_0,{\bf y}_0)\geq 0$,
$\grad c_i^\epsilon (t_0,{\bf y}_0)= 0$
and $\Delta c_i^\epsilon (t_0,{\bf y}_0)\geq 0$ yields a contradiction since $f_i (t_0,c_i^\epsilon (t_0,{\bf y}_0))\geq 0$.
Here, because of the specific boundary condition and
the fact that $\Omega$ has a smooth boundary, the same argument works also if ${\bf y}_0$ is a boundary point.
In the limit $\epsilon \to 0+$ we obtain a nonnegative solution for $\epsilon =0$, hence a
nonnegative solution of the original problem. This finishes the proof since strong solutions
are unique.

Note that non-negativity of the concentrations directly implies $L^\infty$-bounds in the considered isobaric case
due to $0\leq c_i\leq c_{\rm tot}\equiv c_{\rm tot}^0$, which is an important first step for global existence.
\\[2ex]
The considerations in Section~\ref{wellposed} are helpful to verify that the
Maxwell-Stefan multicomponent diffusion is consistent with the second law from thermodynamics.
Indeed, \eqref{flux-repr1} directly yields
\vspace{-0.05in}
\[
- \left[ {\bf J}_i \right]:\left[
\nabla \mu_i
\right]
=
\frac{1}{R\, T}
\Big( (-A_{S|\hat{E}})^{-1} X^{\frac 1 2}
\left[
\nabla \mu_i
\right] \Big)
: \Big( X^{\frac 1 2}
\left[
\nabla \mu_i
\right]\Big) \geq 0,
\]
i.e.\ the entropy inequality is satisfied.
The latter is already well-known in the engineering literature, but with a different
representation of the dissipative term using the individual velocities; cf.\ \cite{Standart}.

For sufficiently regular solutions and under appropriate boundary conditions the entropy inequality can be used as follows.
Let $V({\bf x})=\int_\Omega G({\bf x})\, d{\bf x}$ with $G$ the Gibbs free energy density. Let
\[
W({\bf x},\nabla {\bf x})=- \int_\Omega \left[ {\bf J}_i \right]:\left[ \nabla \mu_i
\right] d{\bf x} \geq 0.
\]
%\\[0.5ex]
Then $(V,W)$ is a Lyapunov couple, i.e.
\[
V({\bf x}(t)) + \int_0^t W({\bf x}(s),\nabla {\bf x}(s)) \, ds \leq V({\bf x}(0)) \mbox{ for } t>0
\]
and all sufficiently regular  solutions.
For ideal systems this yields a priori bounds on the quantities $|\nabla c_i|^2 /c_i$, hence, equivalently, $L_2$-bounds on
$\nabla \sqrt{c_i}$. This type of a priori estimates is well-known in the theory of reaction-diffusion
systems without cross-diffusion; see \cite{Bo15}, \cite{BoPie} and the references given there for more details.\\[2ex]
In the present paper we considered the isobaric and isothermal case because it allows to neglect convective transport and, hence, provides a good starting point. The general case of a multicomponent flow is much more complicated, even in the isothermal case. This case
leads to a Navier-Stokes-Maxwell-Stefan system which will be studied in future work.\\[2ex]
{\bf Acknowledgement.} The author would like to express his thanks to Jan Prüss
(Halle-Wittenberg) for helpful discussions.

\end{document}